\documentclass[12pt]{article}
\usepackage{amsfonts,amssymb}
\usepackage{amsmath}
\usepackage{amsthm}

\newtheorem{lemma}{Lemma}
\newtheorem{prop}{Proposition}
\newtheorem{corollary}{Corollary}

\newtheorem*{propos}{Proposition}
\title{Quadratic Cohomology}
\author{A.~A.~Agrachev\thanks{SISSA, Trieste \& Steklov Math. Inst.,
Moscow}}
\date{}
\begin{document}
\maketitle

\begin{abstract}
 We study homological invariants of smooth families of real quadratic forms as a step towards a
``Lagrange multipliers rule in the large'' that intends to describe topology of smooth maps
in terms of scalar Lagrange functions.
\end{abstract}

\section{Introduction}

Morse theory connects homology of Lebesgue sets and level sets of smooth real functions with critical points of
the functions. The theory is based on a simple observation that a continuous deformation of the function does not
influence the homotopy type of the level and Lebesgue sets for a prescribed value of the function as long as
the value is not critical. Moreover, homology of the Lebesgue set is easier to control than one of the level set.

The same observation holds for level sets of smooth vector-functions. A natural generalization of a Lebesgue
set is the space of solutions of a system of inequalities. The study of systems of inequalities and equations
is partially reduced to the real functions case by the Lagrange multipliers rule. The Lagrange function of a
vector-function $(\phi^1,\ldots,\phi^k)$ is a linear combination \linebreak $p_1\phi^1+\cdots+p_k\phi^k,\ \sum\limits_{i=1}^kp_i^2=1$, where the
coefficients $p_1,\ldots,p_k$ of the linear combination are treated as extra variables, the {\it Lagrange
multipliers}. The vector $0\in\mathbb R^k$ is a critical value of $(\phi^1,\ldots,\phi^k)$ if and only if $0\in\mathbb R$ is a
critical value of the Lagrange function.

The title of the famous Marston Morse's book \cite{Mi} is ``The calculus of variations in the large". This paper is a step towards
a Lagrange multipliers rule in the large. Our first observation, a starting point of the whole story, is
that linearity with respect to the Lagrange multipliers is not important. More precisely, if two Lagrange
functions
$$
f_0(p_1,\ldots,p_k,x)=\sum\limits_{i=1}^kp_i\phi^i_0(x),\quad f_1(p_1,\ldots,p_k,x)=\sum\limits_{i=1}^kp_i\phi^i_1(x)
$$
are connected
by a homotopy $f_t,\ t\in[0,1]$, where $f_t$ are just smooth, not necessary linear with respect to
the Lagrange multipliers and 0 is not a critical value of $f_t$ for all $t\in[0,1]$, then
zero level sets of the vector functions $(\phi^1_0,\ldots,\phi^k_0)$ and $(\phi^1_1,\ldots,\phi^k_1)$
have equal homologies.

A similar property is valid for systems of inequalities; in this case Lagrange multipliers are taken from
the intersection of the sphere with a convex cone. One inequality (like in Morse theory) corresponds
to a point of the sphere. Actually, any point of the sphere of Lagrange multipliers represents a real
function. We can think on usual homology of the space of solutions to the inequality as a kind of generalized
cohomology of the point (different points may have different generalized cohomologies!). Similarly, the generalized
cohomology of a convex subset of the sphere is the usual homology of the space of solutions to the corresponding
system of inequalities. It is easy to extend the construction to more general subsets of the sphere like
submanifolds with boundaries and corners. For the generalized cohomology to have good properties we impose some
regularity conditions. In particular, not all convex subsets of the sphere are available but only those
corresponding to regular systems of inequalities.

The generalized cohomology satisfies a natural modification of the Eilen\-berg--Steen\-rod axioms \cite{EiSt}. The most important
``homotopy axiom" is based on the above property of the homologies of level sets when
regular homotopies of the Lagrange functions are considered.

Such a cohomology theory is determined by the space of function \linebreak $span\{\phi^1,\ldots,\phi^k\}$; different spaces
of functions give different generalized cohomologies. Moreover, as soon as a space of functions and the axioms are fixed we may try to find other cohomology theory that satisfies the same axioms but may be easier to compute.
Such a theory should anyway have an intimate relation to the systems of inequalities and equations. The axioms
imply that the cohomology of a point equals usual homology of space of solutions to the correspondent inequality;
moreover, the cohomology of a convex set vanishes if the correspondent system of inequalities has no solutions.

This general setting is described in Sections 2--4 of the paper. The main results are presented in Sections 5,\,6, where
we build a cohomology theory that satisfies all the axioms in the case the space of functions is the space
of \emph{quadratic forms}. To compute the cohomology we define a spectral sequence $E^r$ (see Section 5) with clear
explicit expressions for all the differentials. The homotopy invariance is proved in Section 6; the proof is
based on the results of \cite{Ag11}.

The page $E^2$ and the differential $d_2$ of the spectral sequence $E^r$ are equal to the page $F^2$ and the
differential $d_2$ of the spectral sequence $F^r$ described in \cite{AgLe}. The sequence $F^r$ converges to the
homology of the space of solutions to the system of quadratic inequalities. We do not know higher differentials
of the sequence $F^r$ and, for the moment, we do not see a reason for two spectral sequences to be equal.
Anyway, this question remains open.

A couple of words on the differentials $d_r$ of the spectral sequence $E^r$. Recall that we deal with families of quadratic forms, i.e. symmetric matrices. Let $\lambda_1(p)\ge\cdots\ge\lambda_n(p)$ be the eigenvalues of the matrix corresponding to the value $p$ of the parameter. A key role in the construction of the differentials is played by the cycles defined by the equations $\lambda_i(p)=\lambda_{i+1}(p)$ in the space of parameters. All differentials $d_r$ are some Massey operations involving these cycles, they are described in Section~5.

The following example shows a flavor of the developed theory and, in particular, the geometric meaning of the
differential $d_3$.
Let us consider the 3-dimensional space $i\mathrm{su}(2)$ of Hermitian $2\times 2-$matrices with zero trace.
An Hermitian $2\times 2-$matrix can be treated as a symmetric real $4\times 4-$matrix commuting with the
multiplication of the vectors in $\mathbb C^2=\mathbb R^4$ by the imaginary unit $i$. Thus $i\mathrm{su}(2)\subset Sym(\mathbb R^4)$, where $Sym(\mathbb R^4)$ is a 10-dimensional space of real symmetric $4\times 4-$matrices.
Given a matrix $S\in Sym(\mathbb R^4)$, let
$\lambda_1(S)\ge\lambda_2(S)\ge\lambda_3(S)\ge\lambda_4(S)$
be its eigenvalues. If $S\in i\mathrm{su}(2)$, then $\lambda_1(S)=\lambda_2(S)=-\lambda_3(S)=-\lambda_4(S)$,
i.\,e. the eigenvalues are double (the eigenspaces are complex lines). Recall that, in general, for an eigenvalue
to be double is a codimension 2 property in $Sym(\mathbb R^4)$

Now take $S_0\in Sym(\mathbb R^4)$ and translate the subspace $i\mathrm{su}(2)$ by $S_0$. We obtain an affine
subspace $S_0+i\mathrm{su}(2) \subset Sym(\mathbb R^4)$. Matrices from this affine subspace are not forced to be Hermitian and the eigenvalues are not necessary double. We set:
$$
C_j^{S_0}=\{H\in i\mathrm{su}(2): \lambda_j(S_0+H)=\lambda_{j+1}(S_0+H)\},\quad j=1,2,3.
$$
For generic $S_0$, $C_j^{S_0}$ are smooth real algebraic curves in the 3-dimensional space $i\mathrm{su}(2)$.

\begin{propos} $C_j^{S_0},\ j=1,2,3,$ are not empty. Moreover, for generic $S_0$, the curve $C_2^{S_0}$ has odd linking numbers with $C_1^{S_0}$ and with $C_3^{S_0}$.
\end{propos}
This proposition is proved in Section 7.

\medskip\noindent {\sl Acknowledgment.} I am grateful to Antonio Lerario for stimulating discussions.

\section{Regular Homotopy}

Let $M$ be a smooth compact manifold. Given $\phi^0,\phi^1\ldots,\phi^k\in C^1(M)$, the system of equations
$\phi^0(x)=\cdots=\phi^k(x)=0$ is {\it regular} if 0 is not a critical value of the map
$$
\varphi=(\phi^0,\ldots,\phi^k)^T:M\to\mathbb R^{k+1}.
$$
A homotopy $\varphi_t=(\phi^0_t,\ldots,\phi^k_t)^T$ is an {\it isotopy} of the system of equations
$\phi^0_t=\cdots=\phi^k_t=0$ if 0 is not a critical value of $\varphi_t,\ \forall\,t\in[0,1]$.

According to the standard Thom lemma, for any isotopy $\varphi_t$ there exists a family of diffeomorphisms
$\Phi_t:M\to M,\ \Phi_0=id$, such that
$$
\varphi^{-1}_t(0)=\Phi_t\left(\varphi^{-1}_0(0)\right),\ \forall\,t\in[0,1].
$$
This is why one uses the term ``isotopy''.
In particular, $\varphi^{-1}_1(0)\cong \varphi_0^{-1},$ \linebreak
$M\setminus \varphi_1^{-1}\cong M\setminus \varphi_0^{-1}$.

Now consider the function $\varphi^*:S^k\times M\to\mathbb R$ defined by the formula
$\varphi^*(p,x)=\langle p,\varphi(x)\rangle$, where $p\in S^k=\{p\in\mathbb R^{k+1}: |p|=1\}$. It is easy to see that 0 is a critical value of $\varphi$ if and only if it is a critical value of $\varphi^*$.

Nothing prevents us from taking any function $f\in C^1(S^k\times M)$. We say that $f$ is regular if 0 is not a critical value of $f$. A homotopy $f_t,\ t\in[0,1]$, such that all $f_t$ are regular we call a {\it regular homotopy}.
We have much more regular homotopies than isotopies. Nevertheless regular homotopy preserves an important information on
the space of solutions to the system of equations.

\begin{prop} Assume that $f_t$ is a regular homotopy and $f_0=\varphi_0^*,\ f_1=\varphi_1^*$. Then
$M\setminus\varphi_0^{-1}(0)$ is homotopy equivalent to $M\setminus\varphi_1^{-1}(0)$.
\end{prop}
\begin{proof}  We set
$$
B_t=\left\{(p,x)\in S^k\times M : f_t(p,x)>0\right\}.
$$
Note that the projections $(p,x)\mapsto x$ restricted to $B_0$ and $B_1$ are fiber bundles over
$M\setminus\varphi_0^{-1}(0)$ and $M\setminus\varphi_1^{-1}(0)$ whose fibers are hemispheres. In particular,
$B_0$ is homotopy equivalent to $M\setminus\varphi_0^{-1}(0)$ and $B_1$ is homotopy equivalent to
$M\setminus\varphi_1^{-1}(0)$.

We need the following Lemma.

\begin{lemma} There exists a smooth family of diffeomorphisms $F_t:S^k\times M\to S^k\times M$ such that
$F_0=id, F_t(B_0)\subset B_t,\ \forall\,t\in[0,1].$
\end{lemma}
\begin{proof}We set $z=(p,x)\in S^k\times M$ and look for a nonautonomous vector field $Z_t(z)$ such that the flow $F_t$ generated by the differential equation $\dot z=Z_t(z)$ has the desired property. It is sufficient to find a field $Z_t$ such that the equality $f_t(z)=0$ implies $\langle d_zf_t,Z_t(z)\rangle>0$.
Moreover, it is sufficient to find such a field locally and then glue local pieces together by a partition
of unity. It remains to mention that we can easily do it locally since 0 is not a critical value of $f_t.$
\end{proof}

Lemma 1 implies that $B_0$ and $B_1$ are homotopy equivalent. Indeed, we can make a time substitution $t\mapsto 1-t$
and find a flow $G_t:S^k\times M\to S^k\times M$ such that  $G_t(B_1)\subset B_{1-t},\ G_0=id$. The maps
$G_1\circ F_1:B_0\to B_0$ and $F_1\circ G_1:B_1\to B_1$ are obviously homotopic to the identity.
\end{proof}

\medskip
Now I would like to extend the just described construction to systems of inequalities. As we'll see very soon,
inequalities are very useful and helpful even if we are mainly interested in the equations. Let $K\subset\mathbb R^{k+1}$
be a closed convex cone. A system of inequalities is a relation $\varphi(x)\in K,\ x\in M$, were, as before,
$\varphi=(\phi^0,\ldots,\phi^k)^T$.
We say that the system of inequalities is {\it regular} (in the strong sense) if
$\mathrm{im}D_x\varphi+K=\mathbb R^{k+1},\ \forall\,x\in\varphi^{-1}(K)$.

We take the dual cone $K^\circ=\{p\in\mathbb R^{k+1}: \langle p,y\rangle\le 0,\ \forall\,y\in K\}$ and
consider the ``manifold with a convex boundary'' $(K^\circ\cap S^k)\times M$.
We say that a subset $V$ of a smooth manifold is a {\it manifold with a convex boundary} if $V$ is covered by coordinate neigborhoods whose intersections with $V$ are diffeomorphic to closed convex subsets of the Euclidean space.
Smooth functions on the manifold with a convex boundary are restrictions of smooth functions on the ambient manifold.
The tangent cone $T_vV$ is the closure of the set of velocities at $v$ of
smooth curves starting from $v$ and contained in $V$.

Let $f:V\to\mathbb R$ be a $C^1$ function. We say that $v\in V$ is a critical point of $f$ if
$\langle d_vf,\xi\rangle\le 0,\ \forall\,\xi\in T_vV$.

\begin{lemma} If the system of inequalities $\varphi(x)\in K$ is regular (in the strong sense), then 0 is not a critical point of $\varphi^*\bigr|_{(K^\circ\cap S^k)\times M}$.
 \end{lemma}
The proof is a straightforward check based on the duality $K^{\circ\circ}=K$; we leave it to the reader.
The inverse statement is not true mainly due to the fact that $T_yK$ is, in general, bigger than $K$.

\medskip
Definitions of regular functions on a manifold with a convex boundary and of regular homotopy for such
functions are verbatim repetitions of the definitions for a manifold without boundary. An obvious modification
of the proof of Proposition~1 gives:
\begin{prop} Assume that $f_t:(K^\circ\cap S^k)\times M\to\mathbb R,\ t\in [0,1]$, is a regular homotopy and
 $f_0=\varphi^*_0\bigr|_{(K^\circ\cap S^k)\times M},\ f_1=\varphi^*_1\bigr|_{(K^\circ\cap S^k)\times M}$. Then
$M\setminus\varphi_0^{-1}(K)$ is homotopy equivalent to $M\setminus\varphi_1^{-1}(K)$.
\end{prop}

{\bf Remark.} Actually, the (obviously modified) proof of Proposition~1 gives more; namely, under conditions
of Proposition~2 the inclusion
$$
(t,B_t)\hookrightarrow\bigcup\limits_{\tau\in[0,1]}(\tau,B_\tau)
$$
of the subspaces of $[0,1]\times V\times M$
is a homotopy equivalence, $\forall\,t\in[0,1]$.

\medskip
 So the homotopy type
of the complement to the space of solutions to the system of inequalities is preserved by regular homotopies.
It happens that homology of the space of solutions is preserved as well.
\begin{prop}  Assume that $f_t:(K^\circ\cap S^k)\times M\to\mathbb R,\ t\in [0,1]$, is a regular homotopy and
 $f_0=\varphi^*_0\bigr|_{(K^\circ\cap S^k)\times M},\ f_1=\varphi^*_1\bigr|_{(K^\circ\cap S^k)\times M}$. Then
the homology groups of $\varphi_0^{-1}(K)$ and $\varphi_1^{-1}(K)$ with coefficients in a field are isomorphic.
\end{prop}
\begin{proof} We start from the case $K \ne -K$, i.\,e. $K$ is not a subspace and the system of inequalities is not
just a system of equations. In this case, $K^\circ\cap S^k$ is contractible and we have the following series
of homotopy equivalences of the pairs:
$$
\left(M,M\setminus\varphi^{-1}_0(K)\right)\sim\left((K^\circ\cap S^k)\times M,B_0\right)\sim
$$
$$
\left((K^\circ\cap S^k)\times M,B_1\right)\sim \left(M,M\setminus\varphi^{-1}_1(K)\right),
$$
where $B_t=\left\{(p,x)\in(K^\circ\cap S^k)\times M: f_t(p,x)>0\right\}$ (see the proof of Proposition~1).
Hence  $H^*(M,M\setminus\varphi^{-1}_0(K))\cong H^*(M,M\setminus\varphi^{-1}_1(K))$. Alexander--Pontryagin
duality completes the proof for this case.

The case of a system of equations is easily reduced to the case just studied  if we add the tautological
inequality $1\ge 0$ to the system. Let us explain this in more detail. If $K$ is a subspace, then we may assume without lack of generality that $K=0$. Now extend the function $f_t$ to $\mathbb R^{k+1}\times M$ as a
degree one homogeneous function with respect to the variable $p$ (keeping the symbol $f_t$ for the extension)
and consider the functions
$$
\bar f_t:(p,\nu,x)\mapsto f_t(p,x)+\nu, \quad |p|^2+\nu^2=1,\ \nu\le 0.
$$
It is easy to see that $\bar f_t$ are regular. To be absolutely rigorous, we have to smooth out $f_t$
at the points $(0,x)$ but, in fact, nothing depends on the way we do it because $\bar f_t$ is far from 0 at
these points.
\end{proof}

\section{Localization}

Let $V$ be a manifold with a convex boundary and $f:V\times M\to\mathbb R$ a $C^1$-function. In this section,
we assume that $M$ is a real-analytic manifold and $f(v,\cdot)$ is a subanalytic function, $\forall\,v\in V$.
It is convenient to think about $f$ as a family of subanalytic functions $f(v,\cdot)$ on $M$ which depends on the
parameter $v\in V$, and we introduce the notation $f_v\doteq f(v,\cdot)$. ``Localization'' in this section
is the localization with respect to the parameter $v$; the variable $x\in M$ remains global.

We say that the family $f_v$ is regular at $v_0\in V$ if the set $\{v_0\}\times f^{-1}_{v_0}(0)$ does not
contain critical points of $f$.
\begin{prop} Assume that the family $f_v,\ v\in V$, is regular at $v_0\in V$. Then $v_0$ has a compact neighborhood
$O_{v_0}$ and centered at $v_0$ local coordinates $\Phi$ such that $U_0\doteq\Phi(O_{v_0})$ is convex and
the function $(f\circ\Phi+t)\bigr|_{\varepsilon U_0\times M}$ is regular for any sufficiently small
nonnegative constants $t,\varepsilon$ one of which is strictly positive.
\end{prop}
\begin{proof} We may assume that $v_0=0$ is the origin of a Euclidean space and $\Phi=id$. Given $a\in C^1(M),\ y\in\mathbb R,$ we set $C_a(y)=\{x\in a^{-1}(y): d_xa=0\}$ If $0\in\mathbb R$ is not a critical value of
$f_0$, i.\,e. $C_{f_0}=\emptyset$, then the statement is obvious; otherwise, for any $x\in C_{f_0}$ there exists
$\nu_x\in U_0$ such that $\left\langle\frac{\partial f}{\partial v}(0,x),\nu_x\right\rangle\ge\alpha>0$, where $\alpha$
is a positive constant. Then, by the continuity, there exists $\delta>0$ such that for any $\tau\in[-\delta,\delta],\
v\in\delta U_0,\ x\in C_{f_v}(\tau)$, there exists $\hat x\in C_{f_0}(0)$ such that
$$
 \left\langle\frac{\partial f}{\partial v}(v,x),\nu_{\hat x}\right\rangle\ge\delta>0 \eqno (1)
$$

Now let $v\in\varepsilon U_0,\ t\in [0,\delta]$, and $x\in C_{f_v}(-t)$; then $d_xf_v=0$ and $|d_xf_0|\le c\varepsilon$ for some constant $c$. We have:
$$
f(0,x)=f(v,x)-\left\langle\frac{\partial f}{\partial v}(v,x),v\right\rangle+o(\varepsilon),
$$
where $\frac{o(\varepsilon)}\varepsilon\to 0$ as $\varepsilon\to 0$ uniformly for all $v\in\varepsilon U_0,\ x\in C_{f_0}(-t),\ t\in[0,\delta]$. Then
$$
\left\langle\frac{\partial f}{\partial v}(v,x),-v\right\rangle= t+f_0(x)-o(\varepsilon)
$$
and, according to (1),
$$
\left\langle\frac{\partial f}{\partial v}(v,x),\varepsilon\nu_{\hat x}-v\right\rangle\ge t+f_0(x)+\varepsilon\delta-o(\varepsilon).
$$
The Lojasevic inequality  \cite{Ku} gives:
$$
|f_0(x)|\le c'|d_xf_0|^{1+\rho}\le c'c^{1+\rho}\varepsilon^{1+\rho},
$$
where $c',\rho$ are positive constants. Hence $\left\langle\frac{\partial f}{\partial v}(v,x),\varepsilon\nu_{\hat x}-v\right\rangle>0$ if $\varepsilon$ is sufficiently small.
\end{proof}

\begin{corollary} Let $V$ be a compact convex set, $0\in V$. Assume that the family $f_v,\ v\in V,$ is regular
 at 0. Then for any sufficiently small $\varepsilon>0$ the homotopy
$$
(t,v,x)\mapsto f(tv,x)+(1-t)\varepsilon,\quad t\in[0,1],\ v\in\varepsilon V,\ x\in M.
$$
between $f\bigr|_{\varepsilon V\times M}$ and the constant family $(v,x)\mapsto f(0,x)+\varepsilon$ is regular.
\end{corollary}

\section{A Cohomology Theory}

Let $M$ be a real-analytic manifold and $\mathcal A\subset C^1(M)$ a set of subanalytic functions. Let
$W\subset V$ be a pair of manifolds with convex boundaries and $f:V\times M\to\mathbb R$ a regular function such that
$f_v\in\mathcal A,\ \forall\,v\in V,$ and $f\bigr|_{W\times M}$ is also regular.

We set $B_f=\{(v,x): v\in V,\ f(v,x)>0\}$ and define
$$
H^\cdot_{\mathcal A}(f_V,f_W)\doteq H^\cdot\left(V\times M,(W\times M)\cup B_f\right),\quad
H^\cdot_{\mathcal A}(f)\doteq H^\cdot_{\mathcal A}(f_V.f_\emptyset).
$$

The pairs of regular functions $(f_V,f_W)$ form a category $\mathfrak F_{\mathcal A}$ with morphisms $\varphi^*:(f^0_{V_0},f^0_{W_0})\mapsto(f^1_{V_1},f^1_{W_1})$,
where $\varphi:V_1\to V_0$ is a $C^1$-map such that $\varphi(W_1)\subset W_0$ and
$f_v^1=f^0_{\varphi(v)},\ \forall\,v\in V_1$. Then $H^\cdot_{\mathcal A}$ is a functor from this category to the
category of commutative groups.

This is a kind of cohomology functor which satisfies natural
modifications of the Steenrod--Eilenberg axioms except for the
dimension axiom. The exactness and excision are obvious and we do
not repeat them. Homotopy axiom deals with $f:[0,1]\times V\times
M\to\mathbb R$ such that $f_{\{t\}\times V}\in\mathfrak
F_{\mathcal A},\forall\,t\in[0,1],$ and claims that the inclusions
$\{t\}\times V\hookrightarrow[0,1]\times V,\ t\in[0,1],$ induce
the isomorphisms of cohomology groups:
$$
H^\cdot_{\mathcal A}\left(f_{[0,1]\times V},f_{[0,1]\times W}\right)\cong
H^\cdot_{\mathcal A}\left(f_{\{t\}\times V},f_{\{t\}\times W})\right).
$$
This simple but not totally obvious fact was explained in Section~2.

The dimension axiom is substituted by the following one: if $V=\{v\}$ is a point then
$$
H^\cdot_{\mathcal A}\left(f_{\{v\}}\right)=H^\cdot\left(M,\{x\in M: f_v(x)>0\}\right).
$$
The ``points'' for us are regular elements of $\mathcal A$ and different points may have different cohomology.

Standard singular cohomology is a special case. Indeed, let the set $\mathcal A$ consist of one point,
$\mathcal A=\{a\}$, and $a(x)<0,\ \forall\,x\in M$. We have:
$$
H^\cdot_{\{a\}}(V,W)=H^\cdot(V,W)\times H^\cdot(M).
$$

Now assume that $\mathcal A+t\subset\mathcal A$ for any nonnegative constant $t$. Given a map $v\mapsto f_v$
from $V$ to $\mathcal A$ we denote by $(f+t)_{[0,c]\times V}$ the map
$(t,v)\mapsto f_v+t,\ t\in[0,c],\ v\in V$. It was proved in Section~3 that for any $v\in V$ there exists a neigborhood $U_v\subset V$ and $\varepsilon>0$ such that the inclusions
$$
U_v\times\{0\}\hookrightarrow U_v\times[0,\varepsilon],\quad
\{v\}\times\{\varepsilon\}\hookrightarrow U_v\times[0,\varepsilon]
$$
induce the isomorphisms of the cohomology groups
$$
H^\cdot_{\mathcal A}\left(f_{U_v}\right)\cong H^\cdot_{\mathcal A}\left((f+t)_{[0,\varepsilon]\times U_v}\right)
\cong H^\cdot_{\mathcal A}\left(f_{\{v\}}+\varepsilon\right).
$$
In other words, cohomology of a ``small neighborhood'' is equal to the cohomology of a ``point''.

Now assume that the cohomology are taken with coefficients in a field and that $\dim M=n$. Then the cohomology
of a ``point''
$$
H^i_{\mathcal A}\left(f_{\{v\}}\right)=H^i\left(M,\{x\in M:
f_v(x)>0\}\right)= H_{n-i}\left(\{x\in M: f_v(x)\le 0\}\right),
$$
$0\le i\le n$, is simply usual homology of the space of solutions to the inequality $f_v(x)\le 0$.

The localization at a point plus the algebraic homology machinery (based on the axioms) gives a good chance to recover the usual
homology of the space of solutions of a system of inequalities from the ones of the individual inequalities of the
form $a(x)\le 0$, where $a\in{\mathcal A}$.
The success is somehow guaranteed if $H^\cdot_{\mathcal A}$ is a unique cohomology theory for $\mathcal A$ that
satisfies the described axioms. On the other hand, any other cohomology theory that satisfies the same axioms
gives additional important invariants of systems of inequalities or equations for functions from $\mathcal A$.

Let me explain it better for regular systems of equations
$$
\phi^0(x)=\cdots=\phi^k(x)=0,\quad \phi^i\in\mathcal A,\ i=0,1,\ldots,k.
$$
An isotopy $\varphi_t=(\phi_t^0,\ldots,\phi_t^k)^T,\ t\in[0,1],$ of such systems is called $\mathcal A$-{\it rigid}
if $\phi^i_t\in\mathcal A$ for all $t\in[0,1]$. In this case, $\varphi_t^*\in\mathfrak F_{\mathcal A}$, where,
recall,
$$
\varphi^*_t:S^k\times M\to\mathbb R,\quad \varphi_t^*(p,x)=\langle p,\varphi_t(x)\rangle.
$$

Let $\hat H_{\mathcal A}$ be a cohomology functor that satisfies our axioms; then, according to the homotopy axiom,
$\hat H_{\mathcal A}(\varphi^*_0)=\hat H_{\mathcal A}(\varphi^*_1)$. In other words, $\hat H_{\mathcal A}$ is an
invariant of the $\mathcal A$-rigid isotopy. Moreover, it is an invariant of regular homotopies in
$\mathfrak F_{\mathcal A}$ that are much more general than  $\mathcal A$-rigid isotopies.

\medskip
 Let $\varphi=(\phi^0,\ldots,\phi^k)^T,\ (\nu,p)\in\mathbb R\times\mathbb R^{k+1},\ x\in M$; we set $\bar\varphi^*(\nu,p,x)=\nu+\langle p,\varphi(x)\rangle$ and denote by $S^{k+1}_-$ the low semi-sphere in
$\mathbb R\times\mathbb R^{k+1},\ S^{k+1}_-=\{(\nu,p): \nu\le 0,\ \nu^2+|p|^2=1\}$.
\begin{prop} If $\varphi^{-1}(0)=\emptyset$, then $\hat H_{\mathcal A}\left(\bar\varphi^*_{S^{k+1}_-}\right)=0$.
\end{prop}
\begin{proof}
Let $c\in\mathbb R,\ B^{k+1}_c=\{(c,p): p\in\mathbb R^{k+1},\ |p|\le 1\}$. Note that $\bar\varphi^*\bigr|_{B_c^{k+1}\times M}$ is a regular function for any $c>0$ (this is true for any smooth map \linebreak $\varphi:M\to\mathbb R^{k+1}$). Moreover,  $\bar\varphi^*_{B_c^{k+1}}$
is regularly homotopic in $\mathfrak F_{\mathcal A}$ to the constant function $c$; indeed, the homothety of the ball $B^{k+1}_c$ to its center along the radii provides us with the desired regular homotopy. Hence $\hat H_{\mathcal A}\left(\bar\varphi^*_{B^{k+1}_c}\right)=0$.

The function $\bar\varphi^*\bigr|_{B_0^{k+1}\times M}$ is regular if and only if $\varphi^{-1}(0)=\emptyset$.
If it is regular, then it is regularly homotopic to $\bar\varphi^*\bigr|_{B_c^{k+1}\times M}$, where $c>0$,
and $H_{\mathcal A}\left(\bar\varphi^*_{B^{k+1}_0}\right)=0$. It remains to note that the homotopy between $\bar\varphi^*_{B^{k+1}_0}$ and $\bar\varphi^*_{S^{k+1}_-}$ induced by the homotopy $(t;\nu,p)\mapsto((1-t)\nu,p),\ t\in[0,1],\ (\nu,p)\in S^{k+1}_-$, is also regular.
\end{proof}

\medskip
Let $M=\mathbb RP^N=\{(x,-x): x\in S^k\}$ and $\mathcal Q(N)$ the space of real quadratic forms on $\mathbb R^{N+1}$
treated as functions on $\mathbb RP^N$. The main goal of this paper is to construct a
cohomology theory $\hat H_{\mathcal Q(N)}$. This is not just an abstract construction:
we give an effective way to compute the cohomology.

In what follows all cochains and cohomologies are with coefficients in $\mathbb Z_2$. We omit the symbol $\mathbb Z_2$
to simplify notations.

\section{A spectral Sequence}

Now we focus on the space $\mathcal Q(N)$ with fixed $N$ and omit the argument $N$ in order to simplify notations.
We denote by the same symbol a quadratic form on $\mathbb R^{N+1}$ and the function
on $\mathbb RP^N$ induced by this form. A quadratic form $q$ induces a regular function on $\mathbb RP^N$ if and only if $\ker q=0$.
More precisely, critical points of $q:\mathbb RP^N\to\mathbb R$ at $q^{-1}(0)$ are exactly $\bar x=(x,-x)\in \mathbb RP^N$
such that $x\in\ker q\cap S^N$.

Some notations. Let $\lambda_1(q)\ge\cdots\ge\lambda_{N+1}(q)$ be the eigenvalues of the symmetric operator associated
to the quadratic form $q\in\mathcal Q$. We set
$$
\Lambda_{j,m}=\{q\in\mathcal Q: \lambda_{j-1}(q)\ne\lambda_j(q)=\lambda_{j+m-1}(q)\ne\lambda_{j+m}(q)\},
$$
$$
\Lambda_{j,m}^0=\{q\in\Lambda_{j,m}: \lambda_j(q)=0\},
$$
$j=1,\ldots,N,\ m=2,\ldots,N-j+2.$
It is well-known that $\Lambda_{j,m}$ is a smooth submanifold of codimension $\frac{m(m+1)}2-1$ in $\mathcal Q$
while $\Lambda_{j,m}^0$ is a codimension 1 submanifold of $\Lambda_{j,m}$ (see \cite[Prop.~1]{Ag11}
for a short proof).

We say that the pair $(f_V,f_W)\in\mathfrak F_{\mathcal Q}$ is in general position if the boundaries $\partial V,\, \partial W$
are smooth and the map $v\mapsto f_v,\ v\in V$, as well as the restrictions of this map to $W,\,\partial V,\,\partial W$ are transversal
to $\Lambda_{j,m}$ and $\Lambda_{j,m}^0$, for $j=1,\ldots,N,\ m=2,\ldots,N-j+2$.

It is sufficient to construct $\hat H(f_V,f_W)$ and check the axioms for the pairs in general position. Indeed, if the the boundaries
$\partial V,\,\partial W$ are smooth, then standard transversality arguments allow to approximate any pair by a pair in the general position.
Moreover, any two sufficiently close approximations are regularly homotopic and have equal cohomology $\hat H$ according to the
homotopy axiom. The cohomology of the given pair is equal, by definition, to the cohomology of a sufficiently close approximation
in general position.

Similar arguments work in the case of nonsmooth boundaries. Given a manifold $V$ with a convex boundary we can always find a smooth vector field transversal to
the boundary $\partial V$. Trajectories of this field passing through $\partial V$ provide us with a tubular neighborhood
of the boundary. Smooth sections of the tubular neighborhood give us smooth approximations of $\partial V$ inside $V$ and we obtain
$\tilde V\subset V$, where $\partial\tilde V$ is a smooth approximation of $\partial V$. The approximation is good if the time to move from
$\partial\tilde V$ to $\partial V$ along trajectories of our transversal vector field is a $C^0$-small semi-concave
function with a uniformly bounded differential (recall that the differential is defined almost everywhere).

It is easy to see that $(f_{\tilde V},f_{\tilde W})\in\mathfrak F_{\mathcal Q}$ for any sufficiently good approximation $\tilde W\subset W,\
\tilde W\subset\tilde V\subset V$. Moreover, natural diffeomorphisms of different tubular neighborhoods induce 
diffeomorphisms homotopic to the identity of good approximations $(\tilde V,\tilde W)$ and natural isomorphisms of cohomologies $\hat H(f_{\tilde V},f_{\tilde W})$.
The cohomology $\hat H(f_V,f_W)$ is equal, by definition, to $\hat H(f_{\tilde V},f_{\tilde W})$, where $(\tilde V,\tilde W)$ is a
sufficiently good approximation of $(V,W)$ by the pair of manifolds with smooth boundaries.

Let $f:V\to\mathcal Q,\ f\in \mathfrak F_{\mathcal Q}$, be in general position\footnote{for simplicity, we keep the symbol $f$ for the map $v\mapsto f_v$} and
$$
V^j_f=\{v\in V: \lambda_j(f(v))>0\},\quad j=1,\ldots,N+1,
$$
a decreasing filtration of $V$ by open subsets. We equip $V$ with a Riemannian metric and take $\varepsilon>0$ so small that $V^j_f$ and $f^{-1}(\Lambda_{j,m})$ are homotopy retracts of their radius $(\dim V)\varepsilon$ neighborhood, $j=1,\ldots,N+1,\ m=2,\ldots,N-j+2.$

Now consider a smooth singular simplex
$\varsigma:\Delta^i\to V$, where $\Delta^i$ is the standard $i$-dimensional simplex. We say
that $\varsigma$ is adapted to $f$ if the diameter of $\varsigma(\Delta^i)$ is smaller than $\varepsilon$ and  the restriction of  $f\circ\varsigma$ to any face $D$
of $\Delta^i$ satisfies the following properties:

(i) $f\circ\varsigma|_D\pitchfork\Lambda_{j,m}$;

(ii) if $\dim D=4$ and $f\circ\varsigma(D)\cap\Lambda_{j,2}\ne\emptyset$ then $f\circ\varsigma(D)\cap\Lambda_{j+1,2}=f\circ\varsigma(D)\cap\Lambda_{j-1,2}=\emptyset,\ j=1,\ldots,N$.

Manifold $V$ admits a triangulation by adapted simplices. The more delicate property (ii) can be achieved because $\bar\Lambda_{j,2}\cap\bar\Lambda_{j+1,2}=\bar\Lambda_{j,3}$ has codimension 5 in $\mathcal Q$.

We denote by $C_{f,i}(V)$ the space of $i$-dimensional singular chains in $V$ generated by the
adapted singular simplices with coefficients in $\mathbb Z_2$. Let $U$ be an open subset of $V$;
then $C_{f,i}(U)$ is a subspace of $C_{f,i}(V)$ generated by singular simplices with values in $U$ and $C^i_f(V,U)$ is the space of linear forms on $C_{f,i}(V)$ that vanish on
$C_{f,i}(U)$. We obtain a cochain complex
$$
\ldots\to C^{i-1}_f(V,U)\stackrel{\delta}{\longrightarrow}C^i_f(V,U)\stackrel{\delta}{\longrightarrow}C^{i+1}_f(V,U)
\to\ldots, \eqno (2)
$$
where $\delta$ is usual coboundary of singular cochains. The cohomology of the complex (2) coincides with
standard cohomology of the pair $(V,U)$ with coefficients in $\mathbb Z_2$:
$\ker\delta/\mathrm{im}\delta= H^\cdot(V,U)$.

We define cocycles $l^j_f\in C^2_f(V),\ j=1,\ldots,N$ as follows: given a singular simplex
$\varsigma\in C_{f,2}(V),\ l^j_f(\varsigma)$ is the intersection number modulo 2 of $f\circ\varsigma$
and $\Lambda_{j,2}$. We have: $l^j_f\smallsmile l^{j+1}_f=0,\ j=1,\ldots,N-1$. Here $\smallsmile$ is the cup
product of singular cochains.
The maps $\ell_j:\varsigma\to\varsigma\smallsmile l^j_f$ define homomorphisms $\ell_j:C^i_f(V,U)\to C^{i+2}_f(V,U)$. We have $\delta\circ\ell_j=\ell_j\circ\delta,\ \ell_j\circ\ell_{j+1}=0$.

Given $\tau>0$ let $V_f^j(\tau)$ be the radius $\tau$ neighborhood of $V^j_f$. We set:
$$
C^i_j(f)=C^i_f\left(V,V_f^j(i\varepsilon)\right),\quad
C^n(f)=\bigoplus\limits_{i+j=n}C^i_{j+1}(f);
$$
then $\ell_j\left(C^i_{j+1}(f)\right)\subset C^{i+2}_j(f)$. Finally, we define the differential $d:C^n(f)\to C^{n+1}(f)$ by the formula $d\bigr|_{C^i_{j-1}}(f)=\delta+\ell_j$.

The cohomology $\hat H_{\mathcal Q}(f)$ is, by definition, the cohomology of the complex
$$
\ldots\to C^{n-1}(f)\stackrel{d}{\longrightarrow}C^n(f)\stackrel{d}{\longrightarrow}C^{n+1}(f)\to\ldots \eqno (3)
$$

{\bf Remark.} A pedantic reader would say that the cochain groups $C^n(f)$ depend on the small parameter $\varepsilon$. It is not hard to see that the cohomologies of complex (3) for different $\varepsilon$ are naturally isomorphic.

\medskip Consider a filtration of the complex $\bigoplus\limits_{n\ge 0}C^n(f)=
\bigoplus\limits_{n\ge 0}\bigoplus\limits_{i\ge 0}C^i_{n-i+1}(f)$ by subcomplexes $\bigoplus\limits_{n\ge 0}\bigoplus\limits_{i\ge\alpha}C^i_{n-i+1}(f),\ \alpha=0,1,\ldots,\dim V$ and the spectral sequence $E^r_{i,j}$ of this filtration converging to $\hat H_{\mathcal Q}(f)$. We have:
$$
E^1_{i,j}=C^i_{j+1}(f),\quad d_1:C^i_{j+1}(f)\to C^{i+1}_{j+1}(f),\ d_1=\delta.
$$
Hence
$$
E^2_{ij}=H^i(V,V_f^{j+1}),\quad d_2: H^i(V,V_f^{j+1})\to H^{i+2}(V,V_f^j).  \eqno (4)
$$
Moreover, the differential (4) is induced by $\ell_j$ and has a very simple explicit expression. Namely, let
$\bar l^j_f\in H^2\left(V,V\setminus f^{-1}(\bar\Lambda_{j,2})\right)$ be the cohomology class of the cocycle
$l^j_f$. Then $d_2$ is the composition of the map
$$
\bar\ell_j:H^i(V,V^{j+1}_f)\to H^{i+2}\left(V,V_f^{j+1}\cup(V\setminus f^{-1}(\bar\Lambda_{j,2}))\right)
$$
defined by the formula $\bar\ell_j(x)=x\smallsmile\bar l^j_f,\ x\in H^i(V,V^{j+1}_f)$, and the homomorphism
$H^{i+2}\left(V,V_f^{j+1}\cup(V\setminus f^{-1}(\bar\Lambda_{j,2}))\right)\to H^{i+2}(V,V^j_f)$ induced
by the inclusion $V^j_f\subset V_f^{j+1}\cup\left(V\setminus f^{-1}(\bar\Lambda_{j,2})\right)$.

We see that $E^2_{i,j}$ and $d_2$ coincide with the second page $F^2_{i,j}$ and the differential
$d_2:F^2_{i,j}\to F^2_{i+2,j-1}$ of the spectral sequence  converging to $H_{\mathcal Q}(f)$ studied in \cite{AgLe}
(see Theorems 25 and 28 of the cited paper). Hence $E^3_{i,j}=F^3_{i,j}$.
Now we are going to give simple explicit
expressions for all differentials $d_r:E^r_{i,j}\to E^r_{i+r,j-r+1},\ r\ge 3$.

Let $\xi\in C^i_{j+1}(f)=E^1_{i,j}$ be a $\delta$-cocycle such that its cohomology class $\bar\xi\in H^i(V,V_f^{j+1})=E^2_{i,j}$ is a $d_2$-cocycle. Then $\xi\smallsmile l^j_f=\delta\eta$, where $\eta\in C^{i+1}_j(f)$. Moreover, $d_3(\bar\xi)$ is the cohomology class of $\eta\smallsmile l^{j-1}_f$ in $H^{i+3}(V,V^{j-1}_f)$ modulo $d_2$-coboundaries while $l^j_f\smallsmile l^{j-1}_f=0$. Hence $d_3(\bar\xi)$ is the Massey product $\langle\bar\xi,\bar l^j_f,\bar l^{j-1}_f\rangle$ combined with an appropriate inclusion homomorphism (see \cite[Ch.~8]{spec} for the definition and basic properties of Massey products).

Now assume that $\xi$ survives in $E^r_{i,j}$, i.\,e. classes of $\xi$ are cocycles for $d_3,\ldots,d_{r-1}$.
The induction procedure implies that $d_r(\xi)$ is the $r$-fold Massey product
$\langle\bar\xi,\bar l^j_f,\ldots,\bar l^{j-r+2}_f\rangle$ combined with appropriate inclusion homomorphisms.

Indeed, since the class of $\xi$ is $d_{r-1}$-cocycle then, according to the induction assumption,
$\langle\bar\xi,\bar l^j_f,\ldots,\bar l^{j-r+3}_f\rangle\ni\delta\zeta$, where $\zeta\in C^{i+r-2}_{j-r+3}(f)$, and $d_r(\bar\xi)$ is the class of $\zeta\smallsmile l^{j-r+2}_f$.

If $\dim V\le k$, then $E^2_{i,j}=0$ for $i>k$. In particular, if $\dim V=3$ then the last possibly nontrivial differential is $d_3$. This differential has a clear geometric meaning that we are going to describe. Assume that $H_1(V;\mathbb Z_2)=0$ and $\partial V$ is connected or empty  (the three-dimensional sphere and ball are available). Then $H_2(V;\mathbb Z_2)=0$ and the linking number mod\,2 of a 1-dimensional cycle in $V$ with a 1-dimensional cycle in $(V,\partial V)$ are well-defined. We have:
$$
d_3: H^0(V,V^{j+1}_f)\longrightarrow H^3(V,V^{j-1}_f). \eqno (*)
$$
Moreover, ranks of $H^0(V,V^{j+1}_f)$ and $H^3(V,V^{j-1}_f)$ are either one or zero.

If both ranks are equal to one, then $d_3$ sends the generator of $H^0(V,V^{j+1}_f)$ to the generator of $H^3(V,V^{j-1}_f)$
multiplied by the linking number of 1-dimensional cycles $f^{-1}(\Lambda_{j,2})$ and $f^{-1}(\Lambda_{j-1,2})$, according to the direct implementation of the above construction.

\medskip Let $W\subset V$ be such that the pair $(f_V,f_W)\in\mathfrak F_{\mathcal Q}$ is in general position
and $\tilde W\supset W$ be an appropriate tubular neighborhood of $W$ such that the pairs $(W,W^j_f)$ are homotopy retracts of $(\tilde W,\tilde W^j_f)$ and $\hat H^\cdot(f_{\tilde W})$ is naturally isomorphic to
$\hat H^\cdot(f_{\tilde W})$. We define:
$$
C^i_j(f_V,f_W)\doteq C^i_j(f_V)\cap C^i_f(V,\tilde W),\quad C^n(f_V,f_W)=\bigoplus\limits_{i+j=n}C^i_{j+1}(f_V,f_W).
$$
The cohomology $\hat H_{\mathcal Q}(f_V,f_W)$ is, by definition, the cohomology of the complex
$$
\ldots\to C^{n-1}(f_V,f_W)\stackrel{d}{\longrightarrow}C^n(f_V,f_W)\stackrel{d}{\longrightarrow}C^{n+1}(f_V,f_W)\to\ldots.
$$
The excision axiom holds automatically while the obvious exact sequence
$$
0\to C^n(f_V,f_W)\to C^n(f_V)\to C^n(f_{\tilde W})\to 0
$$
implies the long exact sequence
$$
\cdots\to \hat H^n_{\mathcal Q}(f_V)\to \hat H^n_{\mathcal Q}(f_W)\to \hat H^{n+1}_{\mathcal Q}(f_V,f_W)\to
\hat H^{n+1}_{\mathcal Q}(f_V)\to\cdots.
$$
If $V=\{v\}$ is a point, then $\hat H_{\mathcal Q}(f_{\{v\}})=H_{\mathcal Q}(f_{\{v\}})$ since the spectral sequence
$E^r_{i,j}$  degenerates in the page $E^2_{i,j}$ in this case.

The homotopy property is automatic for homotopies in the class of functions in the general position. This property
is not at all trivial for homotopies that include functions not in general position. Moreover, this property
is actually the central point of the whole story; we prove it in the next section.

\medskip {\bf Remark.} To be precise, we have to remind that our cochain spaces depend on a small parameter $\varepsilon$. Of course, we simply take $\varepsilon$ smaller each time it is necessary to guarantee that the final result does not depend on $\varepsilon$.

\section{Surgery}

Let $V$ be a manifold with a convex boundary and $f:V\times\mathbb RP^N\to\mathbb R$ a $C^1$-function such that $f_v\in\mathcal Q,\
\forall v\in V$. The function $f$ is regular if and only if for any $(v,\bar x)\in V\times\mathbb RP^N$ such that
$x\in\ker f_v$ there exists $\xi\in T_vV$ such that $\langle\frac{\partial f}{\partial v}(v,\bar x),\xi\rangle>0$.

 We say that $f$ is {\it strongly regular} if for any $v\in V$ such that $\ker f_v\ne 0$ there exists $\xi\in T_vV$
such that  $\langle\frac{\partial f}{\partial v}(v,\bar x),\xi\rangle>0$ for any $x\in\ker f_v\cap S^N$.

In other words, for the regularity to be strong we ask for $\xi$ in the inequality to be one and the same for all
$x\in\ker f_v\cap S^N$. Here is a typical example of a regular but not strongly regular map:
$$
V=\{q\in\mathcal Q: \mathrm{tr}\,q=0,\ |q|\le 1\},\quad f(q,\bar x)=q(x).  \eqno (5)
$$
Here and below we use the following notations: $\mathrm{tr}\,q$ is the trace of the symmetric operator on
$\mathbb R^{N+1}$ associated to $q$, $\langle q_1,q_2\rangle$ is the trace
of the product of the operators associated to $q_1$ and $q_2$, $|q|=\sqrt{\langle q,q\rangle}$.
Strong regularity is violated at $q=0$.

\begin{lemma} If $f\in\mathfrak F_{\mathcal Q}$ is in general position, then $f$ is strongly regular.
 \end{lemma}
\begin{proof} Let $q\in\mathcal Q$ and $\ker q\ne 0$; then $q\in\Lambda^0_{j,m}$ for some $j,m$. It is easy to see
that $T_q\Lambda^0_{j,m}$ is the kernel of the linear map $q'\mapsto q'\bigr|_{\ker q},\ q'\in\mathcal Q$.
Hence the transversality of the map$v'\mapsto f_{v'},\ v'\in V,$ to $\Lambda_{j,m}^0$ at $v\in V$ is equivalent to the surjectivity of the map
$\xi\mapsto\left\langle\frac{\partial f}{\partial v}(v,\cdot),\xi\right\rangle\bigr|_{\ker f_v},\ \xi\in T_vV$, and implies
the existence of $\xi\in T_vV$ such that the quadratic form $\left\langle\frac{\partial f}{\partial v}(v,\cdot),\xi\right\rangle$ is positive definite on $\ker f_v.$
\end{proof}

\medskip {\bf Remark.} We actually proved more than stated: for $f$ to be strongly regular it is sufficient
that the map $v\mapsto f_v,\ v\in M$ is transversal to submanifolds $\Lambda^0_{j,m}$; transversality to
$\Lambda_{j,m}$ is not necessary.

\medskip
We say that a regular homotopy $f_t,\ t\in[0,1]$, is strongly regular if all $f_t$ are strongly
regular.
Example: take $f$ as in (5), $\alpha\in[0,1)$ and the homotopy $f_t=f+t-\alpha$; then $f_t$ is strongly regular
for all $t$ except of $t=\alpha$. We'll show later that this example is in a sense a universal model of a generic
regular but not strongly regular homotopy.

\begin{lemma} Assume that $f_t\in\mathfrak F_{\mathcal Q},\ f_t:V\times\mathbb RP^N\to\mathbb R,\ t\in[0,1],$
 is a strongly regular homotopy. Then there exists a smooth family of diffeomorphisms\footnote{If $\partial V\ne\emptyset$, then $F_t(V)$ may be a proper subset of $V$.}
$F_t:V\to V$, such that
$F_0=id,\ F_t(V^j_{f_0})\subset V^j_{f_t},\ \forall t\in[0,1],\ j=1,\ldots,N+1$.
\end{lemma}
\begin{proof}
The proof is similar to the proof of Lemma~1. It is sufficient to find a smooth vector field $X_t$ on $V$
such that the equality $\lambda_j({f_t}_v)=0$ implies:
$$
\left\langle\frac{\partial f_t}{\partial v}(v,\bar x),X_t(v)\right\rangle >0,\quad \forall\,x\in\ker {f_t}_v\cap S^N.
\eqno(6)
$$
Indeed, fix $t$ and $v$ and consider a trajectory $v(\tau)$ of the flow generated by the field $X_\tau$ such that
$v(t)=v$. Inequality (6) implies that for any smaller than $t$ and sufficiently close to $t$  number $\tau$ the
quadratic form ${f_{\tau}}_{v(\tau)}$ is negative definite on the linear hull of the eigenvectors of the form
${f_t}_v$ corresponding to the eigenvalues $\lambda_j({f_t}_v),\ldots\lambda_{N+1}({f_t}_v)$. Hence
$\lambda_j({f_\tau}_{v(\tau)})<0$, according to the minimax principle for the eigenvalues of a symmetric operator.
We obtain that any  trajectory started in $V^j_{f_0}$ stays in $V^j_{f_t}$ for all $t\in[0,1]$.

The existence of a desired vector field is guaranteed by the strong regularity assumption.
\end{proof}

Lemma~3 immediately implies the following:

\begin{corollary} Strongly regular homotopies preserve the page $E^2_{i,j}$ of the spectral sequence $E^r_{i,j}$
 described in Section~5.
\end{corollary}

A routine transversality technique gives the following:
\begin{prop} Let $\tilde f_t\in\mathfrak F_{\mathcal Q},\ \tilde f_t:V\times\mathcal Q\to\mathbb R,\ t\in[0,1]$ be a regular
homotopy and $\tilde f_0,\tilde f_1$ are in the general position. Then there exists an arbitrarly $C^0$-close
to $\tilde f_t$ regular homotopy $f_t$ such that $f_0=\tilde f_0,\ f_1=\tilde f_1$;
the function $f_t\in\mathfrak F_{\mathcal Q}$ is not in general position only for a finite number of values of the parameter $t\in(0,1)$, and for any $f_t$ that is not in the general position there exists exactly one point $v_t$ where the map
$v\mapsto{f_t}_v,\ v\in V$, is not transversal to a submanifold $\Lambda_{j,m}$ or $\Lambda^0_{j,m}$
Moreover, if $v_t\in int\,V,\ {f_t}_{v_t}\in\Lambda^0_{j,m}$ and the map $v\mapsto {f_t}_v,\ v\in V$, is not transversal
to $\Lambda^0_{j,m}$ at $v_t$, then the following conditions are satisfied:
\begin{itemize} \item The image of the linear map
$\frac{\partial f_t}{\partial v}(v_t,\cdot)\bigr|_{\ker {f_t}_{v_t}}$ from $T_{v_t}V$ into the space of quadratic forms on $\ker {f_t}_{v_t}$
is a subspace of codimension 1 in the space of quadratic forms and the orthogonal
complement to this subspace is generated by
$\frac{\partial}{\partial\tau}\bigr|_{\tau=t}({f_\tau}_{v_t}|_{\ker{f_t}_{v_t}})$.

\item  $\frac{\partial}{\partial\tau}\bigr|_{\tau=t}({f_\tau}_{v_t}|_{\ker{f_t}_{v_t}})$ is a nondegenerate
quadratic form.

\item The Hessian of the map $v\mapsto {f_t}_v|_{\ker{f_t}_{v_t}} ,\ v\in V$ at $v_t$ is a nondegenerate quadratic form on
the kernel of the map $\frac{\partial f_t}{\partial v}(v_t,\cdot)\bigr|_{\ker {f_t}_{v_t}}$.
\end{itemize}
If $v_t\in\partial V$ and the map $v\mapsto f_v,\ v\in\partial V$, is not transversal to $\Lambda^0_{j,m}$, then
the same conditions are satisfied for ${f_\tau}_{\partial V}$ in place of $f_{\tau}$, and the linear map
$\frac{\partial f_t}{\partial v}(v_t,\cdot)\bigr|_{\ker {f_t}_{v_t}}$ from $span\,T_{v_t}V$ into the space of quadratic forms on $\ker {f_t}_{v_t}$ is surjective.
\end{prop}

We are now ready to state a local version of the homotopy invariance property.

\begin{prop} In the setting of Proposition~6, let $t\in(0,1)$ be such that the map $v\mapsto {f_t}_v,\ v\in V$,
is not in general position. Then there exist a neighborhood $O_{v_t}$ of $v_t$ in $V$ and a neighborhood
$o_t$ of $t$ in $(0,1)$ such that the inclusions $\{\tau\}\times O_{v_t}\hookrightarrow o_t\times O_{v_t},\ \tau\in o_t$, induce isomorphisms
$\hat H_{\mathcal Q}(F_{o_t\times O_{v_t}})\cong\hat H_{\mathcal Q}({f_\tau}_{O_t}),$
where $F_{(\tau,v)}\doteq {f_\tau}_v$.
\end{prop}

The general ``global" homotopy invariance property easily follows from Proposition~7. Indeed, a singularity at $(t,v_t)$ does not influence relative cohomologies for the pairs $([0,1]\times V,\,o_t\times O_{v_t}),\
(V, O_{v_t})$ and the inclusion
$$
\left(\{\tau\}\times V,\{\tau\}\times O_{v_t}\right)\hookrightarrow
\left(o_t\times V,o_t\times O_{v_t}\right)
$$
induces an isomorphism
$\hat H_{\mathcal Q}(F_{o_t\times V},F_{o_t\times O_{v_t}})\cong\hat H_{\mathcal Q}({f_\tau}_V,{f_\tau}_{O_t})$.
The exact sequences of the pairs $(F_{o_t\times V},F_{o_t\times O_{v_t}}),\ ({f_\tau}_V,{f_\tau}_{O_t})$ and the
five lemma imply that the inclusion $\{\tau\}\times V\hookrightarrow o_t\times V$ induces an isomorphism
$\hat H_{\mathcal Q}(F_{o_t\times V})\cong\hat H_{\mathcal Q}(f_\tau)$.

\begin{proof}
First assume that the map $v\mapsto {f_t}_v,\ v\in V,$ is transversal to all submanifolds $\Lambda^0_{j,m}$.
Then $f_t$ is strongly regular (see the Remark after Lemma~3). Hence $\tau\mapsto {f_\tau}_{O_{v_t}},\
\tau\in o_t,$ is a strongly regular homotopy for appropriate neighborhoods $O_{v_t},\,o_t$. Moreover, for any
$\tau_0\in o_t$  the maps
$(\tau,v)\mapsto {f_\tau}_v$ and $(\tau,v)\mapsto {f_{\tau_0}}_v$ on $o_t\times O_{v_t}$ are strongly regular homotopic. Hence $F_{o_t\times O_{v_t}}$ and ${f_{\tau_0}}_{O_{v_t}}$ have equal pages $E^2_{i,j}$.

On the other hand, $F_{o_t\times O_{v_t}}$ is regularly homotopic to a constant family $(\tau,v)\mapsto {f_t}_{v_t}+\varepsilon$ according to the general localization result of Section~3. Moreover, this regular
homotopy is strongly regular in the case under consideration and preserves the page $E^2_{i,j}$. The page
$E^2_{i,j}$ of the constant family has only one nonzero column and the same is true for the families
$F_{o_t\times O_{v_t}}$ and ${f_{\tau_0}}_{O_{v_t}}$. In particular, $E^2_{i,j}=E^{\infty}_{i,j}$ are equal
fot these families.

\medskip It remains to study the case when ${f_t}_{v_t}\in\Lambda^0_{j,m}$ and the map
$v\mapsto {f_t}_v,\ v\in V,$ is not transversal to $\Lambda^0_{j,m}$ at $v_t$. Of course it is sufficient
to prove the isomorphism
$\hat H_{\mathcal Q}(F_{o_t\times O_{v_t}})\cong\hat H_{\mathcal Q}({f_\tau}_{O_{v_\tau}})$ for one particular $\tau$ greater than $t$ and one $\tau$ smaller than $t$.

We denote by $Q_t$ the space of quadratic forms on $\ker f_{v_t},\ Q_t=\mathcal Q(m-1)$.
Given $q\in\mathcal Q$, let $E_q\subset\mathbb R^{N+1}$ be the linear hull of the eigenvectors of $q$
corresponding to the eigenvalues $\lambda_j(q),\ldots,\lambda_{j+m-1}(q)$ and $\pi_q:E_q\to\ker f_{v_t}$
be the restriction to $E_q$ of the orthogonal projector of $\mathbb R^{N+1}$ on $\ker f_{v_t}$. Note that
$E_{f_{v_t}}=\ker f_{v_t}$ and $\pi_{f_{v_t}}=id.$ We work in a small neighborhood of $f_{v_t}$ in $\mathcal Q$
and may assume that $E_q$ is transversal to the orthogonal complement of $\ker f_{v_t}$ and $\pi_q$ is
invertible.

Consider a map $\Phi:q\mapsto q\circ\pi^{-1}_q$ from a neighborhood of $f_{v_t}$ to $Q_t$. It is a rational map
and its differential at the point $f_{v_t}$ sends a form $q$ to $q\bigr|_{\ker f_{v_t}}$. Hence $\Phi$ is a
submersion of a neighborhood of $f_{v_t}$ on a neighborhood of the origin in $Q_t$. Moreover,
$\lambda_i(\Phi(q))=\lambda_{j+i-1}(q),\ i=1,\ldots,m$.

We take a sufficiently small neighborhood $O_{v_t}$ of $v_t$ in $V$, a parameter $\tau\in[0,1]$ close to $t$, and define $g_\tau:O_{v_t}\to\mathbb R$ by the formula: ${g_\tau}_v=\Phi({f_\tau}_v)$.
Then $g_\tau\in\mathfrak F_{Q_t}$
and the following equalities are valid\footnote{for simplicity, we keep symbol $g_\tau$ for the map
$v\mapsto {g_\tau}_v$ as in Section~5.}:
$$
V^i_{g_\tau}=V_{f_\tau}^{i+j-1}\cap O_{v_t},\quad g_\tau^{-1}(\Lambda_{i,k})=f_\tau^{-1}(\Lambda_{i+j-1,k})\cap O_{v_t},
$$
$i=1,\ldots,m-1,\ k=2,\ldots,n-i+1.$ Moreover, $O_{v_t}\subset V^{j-1}_f,\ O_{v_t}\cap V_f^{j+m}=\emptyset$.

It follows that the statement of Proposition~7 for $f_{\tau}\in\mathfrak F_{\mathcal Q}$ is equivalent to the same statement for $g_\tau\in\mathfrak F_{Q_\tau}$.

\medskip We have: ${g_\tau}_{v_t}=0$. The family $G:(\tau,v)\mapsto{g_\tau}_v,\ (\tau,v)\in o_t\times O_{v_t}$ is
in general position and is strongly regular homotopic to a constant family $(\tau,v)\mapsto c,\ c>0,$ if $o_t$
and $O_{v_t}$ are sufficiently small. Hence $\hat H_{q_t}(G_{o_t\times O_{v_t}})=0$.

In what follows, we tacitly substitute $o_t$ and $O_{v_t}$ by smaller neighborhoods each time it is necessary
without changing notations. First we study the case $v_t\in int\,V$ and then explain how the case
$v_t\in\partial V$ is reduced to the previous one.

To go ahead we need convenient coordinates in $O_{v_t}$. We put coordinates on $O_{v_t}$ as the product of two balls,
$O_{v_t}=U\times B=\{(u,q): u\in U,\ q\in B\}$, where $U\subset\ker\frac{\partial g_t(v_t)}{\partial v},\ B\subset\mathrm{im}\,\frac{\partial g_t(v_t)}{\partial v}$, in such a way that $v_t=(0,0)$ in our coordinates and
$$
\frac{\partial g_t(v_t)}{\partial v}:(u,q)\mapsto q,\quad u\in\ker\frac{\partial g_t(v_t)}{\partial v},\
q\in\mathrm{im}\,\frac{\partial g_t(v_t)}{\partial v}.
$$
We also set
$q_0=\frac{\partial g_\tau(v_t)}{\partial\tau}\bigr|_{\tau=t}$. Then $B$ is a ball in the hyperplane
$q_0^\perp\subset Q_t$. Recall that $q_0$ is a nondegenerate quadratic form.
Moreover, we assume that the Hessian of the map $v\mapsto g_{tv}$ at $v_t$ is normalized. This means that
$\ker\frac{\partial g_t(v_t)}{\partial v}=span\,U$ is splitted in two subspaces,
$span\,U=\mathbb R^{i_+}\oplus\mathbb R^{i_-}$, and
$$
\frac{\partial^2g_t(0,0)}{\partial u^2}(u)=2(|u_+|^2-|u_-|^2)q_0,\quad u=(u_+,u_-)\in U,\ u_\pm\in\mathbb R^{i_\pm}.
$$

Now we apply a blow-up procedure with a small parameter $\varepsilon>0$. We set:
$$
\varphi^\varepsilon_s(u,q)=\frac 1{\varepsilon^2}g_{t+\varepsilon^2s}(\varepsilon u,\varepsilon^2q),\quad
|s|\le 1,\ (u,q)\in U\times B.
$$
Note that the multiplication of a quadratic form by a positive number does not change the signs and multiplicities
of the eigenvalues. Hence the spectral sequence $E^r_{i,j}$ for $\varphi^\varepsilon_s$ is equal to one
for $(g_\tau)_{(\varepsilon U)\times(\varepsilon^2B)}$ with $\tau=t+\varepsilon^2s$.  We have:
$$
\varphi^\varepsilon_s(u,q)=q+(|u_+|^2-|u_-|^2+s)q_0+O(\varepsilon).
$$
Now fix parameter $s\ne 0$. If $\varepsilon$ is small enough (how small, depends on $s$), then the function $\varphi^\varepsilon_s$ is homotopic
to $\varphi^0_s$ in the class of functions in the general position.

What remains is to prove that $\hat H_{Q_t}(\varphi^0_s)=0$. The following terminology will be useful:
given $\varphi:V\to Q_t,\ \varphi\in\mathfrak F_{Q_t}$, and a homotopy retraction
$h_\tau:V\to V,\ \tau\in[0,1],$ we say that $h_\tau$ is monotone for
$\varphi$ if $V^J_{\varphi\circ h_\tau}\subset V^j_\varphi,\ j=1,\ldots m,\ \tau\in[0,1].$
The homotopy $\tau\mapsto\varphi\circ h_\tau$ induced by a monotone deformation retraction preserves the page $E^2_{i,j},\,d_2$ of the spectral sequence.

We study separately three cases.

\smallskip\noindent{\bf 1.} The quadratic form $q_0$ is sign-indefinite. In this case $q_0^\perp$ contains a
positive definite form $\hat q$. Moreover, if $s$ is sufficiently small then $\hat q+sq_0$ is a positive definite
form. In this case a deformation retraction $h_\tau(u,q)=\left((1-\tau)^{\frac 12}u,\tau\hat q+(1-\tau)q\right)$ is monotone for
$\varphi^0_s$. Indeed,
$$
\varphi^0_s(h_\tau(u,q))=\tau(\hat q+sq_0)+(1-\tau)\left(q+(|u_+|^2-|u_-|^2+s)q_0\right). \eqno (7)
$$
The signature of a quadratic form (i.\,e. the numbers of positive and negative eigenvalues) does not change under
a linear change of coordinates in $\mathbb R^m$, although the eigenvalues do change. Take coordinates such that
the form $\hat q+sq_0$ is represented by a scalar matrix. In these coordinates, eigenvalues of the form (7)
are linear functions of $\tau$. We have: $\varphi^0_s(h_1(u,q))\equiv\hat q+sq_0$. Hence $E^2_{i,j}=0$.

\smallskip\noindent{\bf 2.} The quadratic form $sq_0$ is positive definite. Then the deformation retraction
 $h_\tau(u,q)=\left((1-\tau)^{\frac 12}u,(1-\tau)q\right)$ is monotone for
$\varphi^0_s$ and $\varphi^0_s(h_1(u,q))\equiv sq_0$. Hence $E^2_{i,j}=0$.

\smallskip\noindent{\bf 3.} The quadratic form $sq_0$ is negative definite. In this case, the page $E^2_{i,j}$ is very far
from being zero. We already mentioned that the transformation of $Q_t$ induced by a linear change of coordinates
in $\mathbb R^m$ does not change the signs of eigenvalues and thus the groups $E^2_{i,j}$ of the spectral
sequences associated to elements of $\mathfrak F_{Q_t}$. It is important that the differentials $d_2$ do not
change as well. The last statement needs a justification since the submanifolds $\Lambda_{j,2}\subset Q_t$ do
depend on the choice of coordinates in $\mathbb R^m$. The differential $d_2$ of the spectral sequence $E^r_{i,j}$
does not depend on the choice of coordinates because it is equal to the differential $d_2$ of the spectral sequence $F^r_{i,j}$ constructed in \cite{AgLe} (see Section~5), and $F^r_{i,j}$ is the Leray spectral sequence of a map that respects
changes of coordinates.

Now take coordinates in $\mathbb R^m$ such that the form $q_0$ is represented by a scalar matrix. Then $B$ is a ball in the space of symmetric matrices with zero trace. If $q_0>0$, then the deformation retraction
$(u_+,u_-,q)\mapsto(u_+,(1-\tau)u_-,q),\ \tau\in[0,1],$ is monotone for $\varphi^0_s$. Similarly, if $q_0<0$,
then the deformation retraction
$(u_+,u_-,q)\mapsto((1-\tau)u_+,u_-,q),\ \tau\in[0,1],$ is monotone.

The next lemma completes the proof of Proposition 7 in the case $v_t\in int\,V$,

\begin{lemma} Let $0<s<1,$
$$
U=\{u\in\mathbb R^k: |u|^2\le 2\},\quad \mathbb B=\{q\in\mathcal Q:\mathrm{tr}\,q=0,\ \|q\|\le 1\},
$$
and the map $\varphi:U\times B\to\mathcal Q,\ \varphi\in\mathfrak F_{\mathcal Q}$, is defined by the formula:
$\varphi(u,q)=q+|u|^2-s$. Then the page $E^3_{i,j}$ of the spectral sequence $E^r_{i,j}$ associated to $\varphi$ is zero.
\end{lemma}
\begin{proof}We have to prove that the cochain complex $(E^2,d_2)$ is exact. It is not at all obvious but it is actually proved in
\cite[Th.\,2]{Ag11}. Indeed, let us show that the complex $(E^2,d_2)$ can be naturally identified with complex (1) from \cite{Ag11}, where $n=N+1$.

We set: $M^j=\{q\in\mathbb B: \|q\|=1,\ \lambda_{N-j+1}(q)\ne\lambda_{N+1}(q)\}$,
like in \cite{Ag11} (note that the eigenvalues have the reversed ordering in \cite{Ag11}).
Recall that $E^2_{i,j}=H^i(V,V_\varphi^{j+1})$, where $V=U\times\mathbb B$. A simple homotopy that moves only eigenvalues of symmetric matrices
keeping fixed the eigenvectors gives a homotopy equivalence of pairs:
$$
\left(U\times\mathbb B,V^{j+1}_\varphi\right)\cong\left(U\times\mathbb B,(U\times M^{N-j})\cup(\partial U\times\mathbb B)\right).
$$
Hence $E^2_{i,j}=H^{i-k}(\mathbb B,M^{N-j})$; moreover, natural isomorphism of $E^2_{\cdot,\cdot}$ and $H^{\cdot-k}(\mathbb B,M^{N-\cdot})$ transforms $d_2$
in the differential of the exact complex (1) from \cite{Ag11}.
\end{proof}

\medskip Let $v_t\in\partial V$; we consider the maps $g_\tau|_{\partial V}$, take appropriate coordinates,
and apply the blow-up procedure as we did for $g_\tau$ in the case of an interior point $v_t$. We arrive to the map
$\varphi^0_s:(u,q)\mapsto q+(|u_+|^2-|u_-|^2+s)q_0$
extended to $U\times B^+$ or $U\times B^-$, where $B^\pm$ is the intersection of a ball in $Q_t$ with the half-space
$\{q\in Q_t:\pm\langle q,sq_0\rangle\ge 0\}$. We denote these extensions by $\varphi^{\pm}_s$.
What remains is to prove that $\hat H_{Q_t}(\varphi^{\pm}_s)=0$.

If $sq_0$ is not negative definite and $|s|$ is sufficiently small, then simple monotone deformation retractions transform $\varphi^{\pm}_s$
into a positive constant. The same is true for $\varphi^+_s$ with a negative definite $sq_0$. The only remaining possibility is
$\varphi^-_s$ with a negative definite $sq_0$. In this case, a deformation retraction
$h_\tau(u,q)=\left(u,q-\tau\frac{\langle q,q_0\rangle}{|q_0|^2}q_0\right),\ \tau\in[0,1]$, is monotone and transforms $\varphi^-_s$ in the already studied $\varphi^0_s$
defined on $U\times B.$
\end{proof}

\medskip {\bf Remark.} We have shown that local disturbance in the page $E^2$ caused by a violation of the strong regularity during a regular homotopy
is totally calmed in the page $E^3$. However, this fact does not imply regular homotopy invariance of $E^3$ because the complexes $E^2,d_2$ do not
satisfy the exact sequence ``axiom" and invariance of their local cohomologies does not imply invariance of the global ones.

\section{An example}

Let $\mathbb H$ be the quaternion algebra, $\mathbb H=\mathbb R\oplus\mathbb R^3$, where $\mathbb R$ is the
real line and $\mathbb R^3$ is the space of purely imaginary quaternions,
$\mathbb R^3=\{x\in\mathbb H: \bar x=-x\}$. We take $a\in\mathbb R^3\setminus\{0\}$ and consider a quadratic map
$\varphi:\mathbb H\to\mathbb R^3$ defined by the formula $\varphi(x)=\bar xax$. Then $|\varphi(x)|=|a||x|^2$.
In particular, $\varphi^{-1}(x)=0$. The restriction of $\varphi$ to $S^3$ is just adjoint representation of the
group $\mathrm{SU}(2)=S^3$ and a realization of the Hopf bundle $S^3\to S^2$. Now consider a family of quadratic forms
$\varphi^*_p\in\mathcal Q(3),\ p\in B^3=\{p\in\mathbb R^3: |p|\le 1\}$, where \linebreak
$\varphi^*_p(x)=\langle p,\varphi(x)\rangle$; then $\varphi^*\in\mathfrak F_{\mathcal Q(3)},\ \hat H_{\mathcal Q(3)}(\varphi^*)=0$.

We have $\mathbb H=\mathbb C\oplus j\mathbb C=\mathbb C^2$. Quadratic forms $\varphi^*_p$ are thus real quadratic
forms on $\mathbb C^2$. It is easy to see that they are Hermitian quadratic forms whose Hermitian matrices have
zero traces. In other words, $span\{\varphi^*_p: p\in B^3\}=i\mathrm{su}(2)$. Eigenspaces of the symmetric operators associated to
$\varphi^*_p$  are complex lines in $\mathbb R^4$; hence the eigenvalues are double and we have
$$
\lambda_1(\varphi^*_p)=\lambda_2(\varphi^*_p)=-\lambda_3(\varphi^*_p)=-\lambda_4(\varphi^*_p),
$$
$$
V^1_{\varphi^*}=V^2_{\varphi^*}=B^3\setminus\{0\},\quad V^3_{\varphi^*}=V^4_{\varphi^*}=\emptyset.
$$
Let $\varsigma$ be a small quadratic form, then $\phi^*-\varsigma$ is regularly homotopic to $\varphi^*$ and
$\hat H_{\mathcal Q(3)}(\varphi^*-\varsigma)=0$. Moreover, $\varphi^*-\varsigma$ is in general position for
almost every $\varsigma$.

Assume that $\varsigma$ is positive definite; then $V^1_{\varphi^*-\varsigma},\ V^2_{\varphi^*-\varsigma}$ are complements
to (small) contractible neighborhoods of 0, $V^3_{\varphi^*-\varsigma}=V^4_{\varphi^*-\varsigma}=\emptyset$.
Indeed, the number of positive eigenvalues of the operator associated to a quadratic form does not depend on the
choice of the Euclidean structure. If we choose a form $\frac 1\varepsilon\varsigma$ as the Euclidean structure, then $\lambda_i(\varphi_p^*-\varsigma)=\lambda_i(\varphi^*_p)-\varepsilon$.

The page $E^2$ of the spectral sequence $E^r$ for $\varphi^*-\varsigma$ has the form:
$$
\begin{array}{|c|c|c|c}
\mathbb Z_2&0&0&0\\
\mathbb Z_{2}&0&0&0\\
0&0&0&\mathbb Z_{2}\\
0&0&0&\mathbb Z_2\\
\hline
\end{array}
$$
Hence the differentials
$
d_3: E^2_{0,j+1}\to E^2_{3,j-1},\ j=2,3,
$
are not zero. We are in the situation described in Section 5 (see the paragraph with formula $(*)$ and the next paragraph). It follows that the linking number mod 2 of $(\varphi^*-\varsigma)^{-1}(\Lambda_{2,2})$ with
$(\varphi^*-\varsigma)^{-1}(\Lambda_{1,2})$ and with $(\varphi^*-\varsigma)^{-1}(\Lambda_{3,2})$ are nonzero.

The Proposition stated in the Introduction can be easily derived from this fact. We start from the case of generic
$S_0$. First of all, $C_i^{S_0+tI}=C_i^{S_0}$ for any scalar matrix $tI$. Hence we may assume that $S_0$ is the matrix of a negative definite quadratic form. It is sufficient to compute linking numbers of $C_2^{S_0}$ with
$C_1^{S_0}$ and with $C_3^{S_0}$ in a very big ball $\frac 1\varepsilon B^3$. Multiplication by $\varepsilon$ transforms $C_j^{S_0}$ into $C_j^{\varepsilon S_0}=(\varphi^*-\varsigma)^{-1}(\Lambda_{j,2}),\
j=1,2,3,$ where $\varsigma$ is the quadratic form represented by the matrix $-\varepsilon S_0$.

We have proved the statement about linking numbers in the case of generic $S_0$. Now take any $S_0$ and present
it as the limit of a sequence of generic ones, $S_0=\lim\limits_{n\to\infty}S_0^n$. Any limiting point of the
sequence of sets $C_j^{S_0^n}$ as $n\to\infty$ belongs to $C_j^{S_0}$. The curves $C_2^{S_0^n}$ are uniformly bounded, hence $C_2^{S_0}\ne\emptyset$. The curves $C_1^{S_0^n}$ and $C_3^{S_0^n}$ are linked with $C_2^{S_0^n}$
and cannot escape to infinity; hence $C_1^{S_0}$ and $C_3^{S_0}$ are also nonempty.

\section{Informal discussion}

The anonymous referee asked me to say more about global features of the Lagrange multipliers even if we do not have yet a general conventional theory. Indeed, Arnold journal encourages informal discussions, and I'll try to do it.

Let $F:U\to M$ be a smooth map from one smooth manifold to another one. Given a critical point $u\in U$ of this map, a Lagrange multiplier is a nonzero covector $\lambda\in T^*_{F(u)}M$, which annihilates the image of the differential $D_xF:T_uU\to T_{F(u)}M$. In other words, Lagrange multipliers are solutions of the equation
$\lambda D_uF=0$ where the pair $(\lambda,u)$ is taken from the total space of the vector bundle
$F^*(T^*M)$ with a removed zero section. The equation is homogeneous on the fibers of the bundle.

The traditional nonhomogeneous ``affine" version of this equation concerns the case $M=\mathbb R\times N,\
F=(\varphi,\Phi)$, where $\phi:U\to\mathbb R$ is treated as a ``functional" and $\Phi:U\to M$ defines constraints.
The Lagrange multiplier is now an element of $T^*_{F(u)}(\mathbb R\times N)=\mathbb R\oplus T^*_{\Phi(u)}N$.
Let $u$ be a regular point of $\Phi$; then $u$ is critical for $F$ if and only if it is a critical point of
$\varphi$ restricted to the level set of $\Phi$. The first (scalar) component of the Lagrange multiplier does not vanish in this case and can be normalized. We set this scalar to be equal to $(-1)$ and obtain the equation:
$\lambda D_u\Phi=d_u\varphi,\ \lambda\in T^*_{\Phi(u)}M$.
The pair $(\lambda,u)$ belongs to $\Phi^*(T^*M)$ and $\lambda$ is also called the Lagrange multiplier.
Both homogeneous and ``affine" versions can be treated similarly.

The map $(\lambda,u)\mapsto\lambda D_uF$ is transversal to the zero section of $F^*(T^*M)$ for generic $F$.
If it is transversal then we say that $F$ is a Morse map. Indeed, for $M=\mathbb R$ this just a usual Morse function. For a Morse map $F$, solutions of the equation $\lambda D_uF=0$ form a smooth $(\dim M)$-dimensional submanifold $C_F$ of $F^*(T^*M)$ (or a $(\dim M-1)$-dimensional submanifold of the projectivization of this vector bundle).

In other words, Lagrange multipliers resolve singularities of the set of critical points. Moreover, the map
$F^c:(\lambda,u)\mapsto\lambda,\ (\lambda,u)\in C_F$
is a Lagrangian immersion of $C_F$ into the manifold $T^*M$ endowed with the standard symplectic structure. Similarly for the affine version, and all that is almost a tautology (see \cite{AgGa} for some details).
I find it wonderful that Lagrange multipliers form a Lagrange submanifold! Both objects are named after Lagrange but they look very different at the first glance.

The idea is to recover interesting homological invariants of $F$ in terms of the Lagrange multipliers sitting in $T^*M$. We would like to develop a theory, which is efficient when $M$ has a modest dimension while $U$ can be huge. The applications most interesting for us concern constrained variational problems where $U$ is an infinite dimensional Hilbert or Banach manifold.

The results of this paper can be easily interpreted as a desired theory for homogeneous quadratic maps. Why do we think that a good theory can be developed in the general setting as well? To any $(\lambda,u)\in C_F$ we
associate the Hessian $\lambda Hess_uF:\ker D_uF\to\mathbb R$ that is a quadratic form on $\ker D_uF$. If $M=\mathbb R$ then critical points of $F$ are isolated, the Hessians of $F$ at these points are nondegenerate quadratic forms and inertia indices of these quadratic forms are crucial local invariants used by the Morse theory to estimate homology of the Lebesgue sets and level sets of $F$. If $\dim M>1$ then critical
points are not isolated and $\lambda Hess_uF$ may be degenerate for some $(\lambda,u)\in C_F$.

There is an important duality between the quadratic form $\lambda Hess_uF$ and the image of the tangent space
$T_{(\lambda,u)}C_F$ under the Lagrangian immersion $F^c:(\lambda,u)\mapsto\lambda$. Let
$J_\lambda=F^c_*(T_{(\lambda,u)}C_F)$ and $\pi:T^*M\to M$ be the canonical projection.
It is easy to check that $\lambda Hess_uF$ is degenerate if and only if $\pi_*|_{J_\lambda}$ is degenerate
and $\dim\ker \lambda Hess_uF=\dim\ker\left(\pi_*|_{J_\lambda}\right)$. Moreover, for any continuous curve
$(\lambda_t,u_t)\in C_F,\ t\in [0,1]$, such that $\lambda_0Hess_{u_0}F$ and $\lambda_1 Hess_{u_1}F$ are nondegenerate, the difference of inertia indices of these quadratic forms is equal to the Arnold--Maslov index
of the curve $t\mapsto J_{\lambda_t}$. In other words, Arnold--Maslov cocycle of the Lagrangian immersion
equals the co-boundary of of the inertia index of the Hessian.

It is natural to expect that homological invariants of the Lagrangian immersion properly glue together the Hessians
corresponding to different points of one and the same connected component of $C_F$ to give such a connected component the role played by the isolated critical point in the usual Morse theory.

The framework is indeed rather similar to one studied in this paper. Let ${\bf L}_\lambda$ be the Lagrange Grassmannian of all Lagrangian subspaces of the symplectic space $T_\lambda(T^*M)$. This Lagrange Grassmannian has a distiguished element $\Pi_\lambda=T_{\lambda}(T^*_{\pi(\lambda)}M)$ (the tangent space to the fiber) and is, actually, a natural compactification of the space of quadratic forms on $\Pi_\lambda$ (see, for instance, \cite{Ar} or \cite{AgGa}). The subspace $J_{\lambda}$ is also an element of ${\bf L}_{\lambda}$.

Given $\Lambda\in{\bf L}_\lambda$, we have: $\ker\pi_*|_\Lambda=\Lambda\cap\Pi_\lambda$. The set of all Lagrangian subspaces which have a nontrivial intersection with $\Pi_\lambda$ is called ``the train of $\Pi_\lambda$". So the Hessian changes its inertia index exactly when $J_{\lambda}$ passes the train. On the other hand, the train is the compactification of the space of degenerate quadratic forms on $\Pi_\lambda$ (see \cite{Ar}). It looks like we always speak about one and the same story... .

\end{document}